\documentclass[11pt,twoside]{article}
\usepackage{amsmath, amsthm, amscd, amsfonts, amssymb, graphicx, color}

\date{}

\begin{document}

\centerline{}

\centerline {\Large{\bf Generalized fusion frame in Quaternionic  Hilbert spaces }}

\newcommand{\mvec}[1]{\mbox{\bfseries\itshape #1}}
\centerline{}
\centerline{\textbf{Prasenjit Ghosh}}
\centerline{Department of Pure Mathematics, University of Calcutta,}
\centerline{35, Ballygunge Circular Road, Kolkata, 700019, West Bengal, India}
\centerline{e-mail: prasenjitpuremath@gmail.com}
\centerline{}

\newtheorem{Theorem}{\quad Theorem}[section]

\newtheorem{definition}[Theorem]{\quad Definition}

\newtheorem{theorem}[Theorem]{\quad Theorem}

\newtheorem{remark}[Theorem]{\quad Remark}

\newtheorem{corollary}[Theorem]{\quad Corollary}

\newtheorem{note}[Theorem]{\quad Note}

\newtheorem{lemma}[Theorem]{\quad Lemma}

\newtheorem{example}[Theorem]{\quad Example}

\newtheorem{result}[Theorem]{\quad Result}
\newtheorem{conclusion}[Theorem]{\quad Conclusion}

\newtheorem{proposition}[Theorem]{\quad Proposition}

\centerline{}

\begin{abstract}
\textbf{\emph{We introduce the notion of a generalized fusion frame in quaternionic Hilbert space.\,A characterization of generalized fusion frame in quaternionic Hilbert space with the help of frame operator is being discussed.\,Finally, we construct \,$g$-fusion frame in quaternionic Hilbert space using invertible bounded right \,$\mathfrak{Q}$-linear operator on quaternionic Hilbert space. }}
\end{abstract}



\section{Introduction and preliminaries}
 
\smallskip\hspace{.6 cm}In 1952, Duffin and Schaeffer \cite{Duffin} introduced frames for Hilbert spaces to study some fundamental problems in non-harmonic Fourier series.\;Later on, after some decades, frame theory was popularized by Daubechies et al. \cite{Daubechies}.\;At present, frame theory has been widely used in signal and image processing, filter bank theory, coding and communications, system modeling and so on.\;Several generalizations of frames  namely, \,$g$-frames \cite{Sun}, fusion frames \cite{Kutyniok}, \,$g$-fusion frames \cite{Ahmadi} etc. have been introduced in recent times.

Sadri et al. \cite{Ahmadi} studied \,$g$-fusion frame in Hilbert space to generalize the theory of fusion frame and \,$g$-frame.\,Let \,$\left\{\,W_{j}\,\right\}_{ j \,\in\, J}$\, be a collection of closed subspaces of a Hilbert space \,$H$\; and \,$\left\{\,v_{j}\,\right\}_{ j \,\in\, J}$\; be a collection of positive weights and for each \,$j \,\in\, J$, \,$\Lambda_{j} \,:\, H \,\to\, H_{\,j}$\, be a bounded linear operator, where \,$J$\; is subset of  integers \,$\mathbb{Z}$.\;Then the family \,$\Lambda \,=\, \{\,\left(\,W_{j},\, \Lambda_{j},\, v_{j}\,\right)\,\}_{j \,\in\, J}$\; is called a generalized fusion frame or a $g$-fusion frame for \,$H$\; respect to \,$\left\{\,H_{j}\,\right\}_{j \,\in\, J}$\; if there exist constants \,$0 \,<\, A \,\leq\, B \,<\, \infty$\; such that
\begin{equation}\label{eq1}
A \;\left \|\,f \,\right \|^{\,2} \,\leq\, \sum\limits_{\,j \,\in\, J}\,v_{j}^{\,2}\, \left\|\,\Lambda_{j}\,P_{\,W_{j}}\,(\,f\,) \,\right\|^{\,2} \,\leq\, B \; \left\|\, f \, \right\|^{\,2}\; \;\forall\; f \,\in\, H,
\end{equation}
where each \,$P_{\,W_{j}}$\, is a orthogonal projection onto the closed subspace \,$W_{\,j}$, for \,$j \,\in\, J$\, and \,$\left\{\,H_{j}\,\right\}_{ j \,\in\, J}$\, is the collection of Hilbert spaces.\,The constants \,$A$\; and \,$B$\; are called the lower and upper bounds of $g$-fusion frame, respectively.\,If \,$A \,=\, B$\; then \,$\Lambda$\; is called tight $g$-fusion frame and if \;$A \,=\, B \,=\, 1$\, then we say \,$\Lambda$\; is a Parseval $g$-fusion frame.\;If \,$\Lambda$\; satisfies only the right inequality of (\ref{eq1}) then it is called a $g$-fusion Bessel sequence with bound \,$B$\; in \,$H$.

Let \,$\Lambda \,=\, \left\{\,\left(\,W_{j},\, \Lambda_{j},\, v_{j}\,\right)\,\right\}_{j \,\in\, J}$\, be a $g$-fusion Bessel sequence in \,$H$\, with a bound \,$B$.\;The synthesis operator \,$T_{\Lambda}$\, of \,$\Lambda$\; is defined as 
\[ T_{\Lambda} \,:\, l^{\,2}\left(\,\left\{\,H_{j}\,\right\}_{ j \,\in\, J}\,\right) \,\to\, H,\]
\[T_{\Lambda}\,\left(\,\left\{\,f_{\,j}\,\right\}_{j \,\in\, J}\,\right) \,=\,  \sum\limits_{\,j \,\in\, J}\, v_{j}\, P_{\,W_{j}}\,\Lambda_{j}^{\,\ast}\,f_{j}\; \;\;\forall\; \{\,f_{j}\,\}_{j \,\in\, J} \,\in\, l^{\,2}\left(\,\left\{\,H_{j}\,\right\}_{ j \,\in\, J}\,\right)\] and the analysis operator is given by 
\[ T_{\Lambda}^{\,\ast} \,:\, H \,\to\, l^{\,2}\left(\,\left\{\,H_{j}\,\right\}_{ j \,\in\, J}\,\right),\; T_{\Lambda}^{\,\ast}\,(\,f\,) \,=\,  \left\{\,v_{j}\,\Lambda_{j}\, P_{\,W_{j}}\,(\,f\,)\,\right\}_{ j \,\in\, J}\; \;\forall\; f \,\in\, H.\]
The $g$-fusion frame operator \,$S_{\Lambda} \,:\, H \,\to\, H$\; is defined as follows:
\[S_{\Lambda}\,(\,f\,) \,=\, T_{\Lambda}\,T_{\Lambda}^{\,\ast}\,(\,f\,) \,=\, \sum\limits_{\,j \,\in\, J}\, v_{j}^{\,2}\; P_{\,W_{j}}\, \Lambda_{j}^{\,\ast}\; \Lambda_{j}\, P_{\,W_{j}}\,(\,f\,)\]and it can be easily verify that 
\[\left<\,S_{\Lambda}\,(\,f\,),\, f\,\right> \,=\, \sum\limits_{\,j \,\in\, J}\, v_{j}^{\,2}\, \left\|\,\Lambda_{j}\, P_{\,W_{j}}\,(\,f\,) \,\right\|^{\,2}\; \;\forall\; f \,\in\, H.\]
Furthermore, if \,$\Lambda$\, is a $g$-fusion frame with bounds \,$A$\, and \,$B$\, then from (\ref{eq1}),
\[\left<\,A\,f,\, f\,\right> \,\leq\, \left<\,S_{\Lambda}\,(\,f\,),\, f\,\right> \,\leq\, \left<\,B\,f \,,\, f\,\right>\; \;\forall\; f \,\in\, H.\]
The operator \,$S_{\Lambda}$\; is bounded, self-adjoint, positive and invertible.\\

In recent times, frames for finite dimensional quaternionic Hilbert spaces were studied by Khokulan et al. \cite{Mk}.\,Sharma and Goel \cite{SKS} introduced frames in a quaternionic Hilbert spaces.\,Various generalization of frame in quaternionic Hilbert space were introduced by S. K. Sharma et al. \cite{AMJ}.

In this paper, we give the notion of a \,$g$-fusion frame in quaternionic Hilbert space and establish a characterization of generalized fusion frame in quaternionic Hilbert space using its frame operator.\,At the end, \,$g$-fusion frames in quaternionic Hilbert spaces using invertible bounded right \,$\mathfrak{Q}$-linear operator on quaternionic Hilbert space are being discussed.

\section{Quaternionic Hilbert space}

We start with this section by giving some basic facts about the algebra of quaternions, right quaternionic Hilbert space and operators on right quaternionic Hilbert spaces.\,The non-commutative field of quaternions \,$\mathfrak{Q}$\, is a four dimensional real algebra with unity.\,In \,$\mathfrak{Q}$, \,$0$\, denotes the null element and \,$1$\, denotes the identity with respect to multiplication.\,It also includes three so-called imaginary units, denoted by \,$i,\, j,\, k$.\,Thus,
\[\mathfrak{Q} \,=\, \left\{\,a_{\,0} \,+\, a_{\,1}\,i \,+\, a_{\,2}\,j \,+\, a_{\,3}\,k \,:\, a_{\,0},\, a_{\,1},\, a_{\,2},\, a_{\,3} \,\in\, \mathbb{R}\,\right\}\]
where \,$i^{\,2} \,=\, j^{\,2} \,=\, k^{\,2} \,=\, -\, 1$; \,$i\,j \,=\, -\, j\,i \,=\, k$; \,$j\,k \,=\, -\, k\,j \,=\, i$\, and \,$k\,i \,=\, -\, i\,k \,=\, j$.\,For each quaternion \,$q \,=\, a_{\,0} \,+\, a_{\,1}\,i \,+\, a_{\,2}\,j \,+\, a_{\,3}\,k \,\in\, \mathfrak{Q}$, the conjugate of \,$q$\, is denoted by \,$\overline{q}$\, and defined by \,$\overline{q} \,=\, a_{\,0} \,-\, a_{\,1}\,i \,-\, a_{\,2}\,j \,-\, a_{\,3}\,k \,\in\, \mathfrak{Q}$.\,Here \,$a_{\,0}$\, is called the real part of \,$q$\, and \,$a_{\,1}\,i \,+\, a_{\,2}\,j \,+\, a_{\,3}\,k$\, is called the imaginary part of \,$q$.\,The modulus of quaternion \,$q$\, is defined as \,$|\,q\,| \,=\, \sqrt{a_{\,0}^{\,2} \,+\, a_{\,1}^{\,2} \,+\, a_{\,2}^{\,2} \,+\, a_{\,3}^{\,2}}$.\,For every non-zero quaternion \,$q \,=\, a_{\,0} \,+\, a_{\,1}\,i \,+\, a_{\,2}\,j \,+\, a_{\,3}\,k \,\in\, \mathfrak{Q}$, there exists a unique inverse \,$q^{\,-\, 1}$\, in \,$\mathfrak{Q}$\, as 
\[q^{\,-\, 1} \,=\, \dfrac{\overline{q}}{|\,q\,|^{\,2}} \,=\, \dfrac{a_{\,0} \,-\, a_{\,1}\,i \,-\, a_{\,2}\,j \,-\, a_{\,3}\,k}{\sqrt{a_{\,0}^{\,2} \,+\, a_{\,1}^{\,2} \,+\, a_{\,2}^{\,2} \,+\, a_{\,3}^{\,2}}}\,.\]

\begin{definition}\cite{RGV}
A right quaternionic vector space \,$\mathbb{H}^{\,R}\,(\,\mathfrak{Q}\,)$\, is a linear vector space under right scalar multiplication over the filed of quaternionic \,$\mathfrak{Q}$, i.\,e.,
\begin{align*}
&\mathbb{H}^{\,R}\,(\,\mathfrak{Q}\,) \,\times\, \mathfrak{Q} \,\to\, \mathbb{H}^{\,R}\,(\,\mathfrak{Q}\,) \,\Rightarrow\, (\,u,\, q\,) \,\to\, u\,q  
\end{align*}
and for each \,$u,\, v \,\in\, \mathbb{H}^{\,R}\,(\,\mathfrak{Q}\,)$\, and \,$p,\, q \,\in\, \mathfrak{Q}$, the right scalar multiplication satisfying the following properties:
\begin{align*}
&(\,u \,+\, v\,)\,q \,=\, u\,q \,+\, v\,q\,,\, \, u\,(\,p \,+\, q\,) \,=\, u\,p \,+\, u\,q\,,\, \,v\,(\,p\,q\,) \,=\, (\,v\,p\,)\,q.
\end{align*}  
\end{definition} 

\begin{definition}\label{1.df.01}\cite{RGV}
A right quaternionic inner product space \,$\mathbb{H}^{\,R}\,(\,\mathfrak{Q}\,)$\, is a right quaternionic vector space together with the binary mapping \,$\left<\,\cdot,\, \cdot\,\right> \,:\, \mathbb{H}^{\,R}\,(\,\mathfrak{Q}\,) \,\times\, \mathbb{H}^{\,R}\,(\,\mathfrak{Q}\,) \,\to\, \mathfrak{Q}$\, which satisfies the following properties:
\begin{itemize}
\item[$(i)$] $\left<\,v,\, v\,\right> \,>\, 0$\, if \,$v \,\neq\, 0$
\item[$(ii)$] $\overline{\left<\,u,\, v\,\right>} \,=\, \left<\,v,\, u\,\right>$\, for all \,$u,\, v \,\in\, \mathbb{H}^{\,R}\,(\,\mathfrak{Q}\,)$.
\item[$(iii)$]$\left<\,u,\, v_{\,1} \,+\, v_{\,2}\,\right> \,=\, \left<\,u,\, v_{\,1}\,\right> \,+\, \left<\,u,\, v_{\,2}\,\right>$\, for all \,$u,\, v_{\,1},\, v_{\,2} \,\in\, \mathbb{H}^{\,R}\,(\,\mathfrak{Q}\,)$.
\item[$(iv)$]$\left<\,u,\, v\,q\,\right> \,=\, \left<\,u,\, v\,\right>\,q$\, for all \,$u,\, v \,\in\, \mathbb{H}^{\,R}\,(\,\mathfrak{Q}\,)$\, and \,$q \,\in\, \mathfrak{Q}$.   
\end{itemize} 
\end{definition}     

Let \,$\mathbb{H}^{\,R}\,(\,\mathfrak{Q}\,)$\, be a right quaternionic inner product space with respect to the right quaternionic inner product \,$\left<\,\cdot,\, \cdot\,\right>$.\,Define the quaternionic norm \,$\left\|\,\cdot\,\right\| \,:\, \mathbb{H}^{\,R}\,(\,\mathfrak{Q}\,) \,\to\, \mathbb{R}^{\,+}$\, on \,$\mathbb{H}^{\,R}\,(\,\mathfrak{Q}\,)$\, by 
\begin{align}
\left\|\,u\,\right\| \,=\, \sqrt{\left<\,u,\, u\,\right>},\, \,u \,\in\, \mathbb{H}^{\,R}\,(\,\mathfrak{Q}\,).\label{1.eq1.01}
\end{align}

\begin{definition}\cite{RGV}
The right quaternionic inner product space \,$\mathbb{H}^{\,R}\,(\,\mathfrak{Q}\,)$\, is called a right quaternionic Hilbert space if it is complete with respect to above norm (\ref{1.eq1.01}).
\end{definition}

\begin{theorem}(Cauchy-Schwarz inequality)\,\cite{RGV}
Let \,$\mathbb{H}^{\,R}\,(\,\mathfrak{Q}\,)$\, be a right quaternionic Hilbert space.\,Then
\[\left|\,\left<\,u,\, v\,\right>\,\right|^{\,2} \,\leq\, \left<\,u,\, u\,\right>\,\left<\,v,\, v\,\right>\,,\;\text{for all}\; \;u,\, v \,\in\, \mathbb{H}^{\,R}\,(\,\mathfrak{Q}\,).\] 
\end{theorem} 
The quaternionic norm defined in (\ref{1.eq1.01}) satsifies the following properties:
\begin{itemize}
\item[$(i)$] $\|\,u\,\| \,=\, 0$\, for some \,$u \,\in\, \mathbb{H}^{\,R}\,(\,\mathfrak{Q}\,)$, then \,$u \,=\, 0$.
\item[$(ii)$] $\|\,u\,q\,\| \,=\, |\,q\,|\,\|\,u\,\|$\, for all \,$u \,\in\, \mathbb{H}^{\,R}\,(\,\mathfrak{Q}\,)$\, and \,$q \,\in\, \mathfrak{Q}$.
\item[$(iii)$]$\|\,u \,+\, v\,\| \,\leq\, \|\,u\,\| \,+\, \|\,v\,\|$\, for all \,$u,\, v \,\in\, \mathbb{H}^{\,R}\,(\,\mathfrak{Q}\,)$.   
\end{itemize}

\begin{definition}\cite{RGV}
Let \,$\mathbb{H}^{\,R}\,(\,\mathfrak{Q}\,)$\, be a right quaternionic Hilbert space and \,$V$\, be a subset of \,$\mathbb{H}^{\,R}\,(\,\mathfrak{Q}\,)$.\,Define
\begin{itemize}
\item[$(i)$]\,$V^{\,\bot} \,=\, \left\{\,v \,\in\, \mathbb{H}^{\,R}\,(\,\mathfrak{Q}\,) \,:\, \left<\,v,\, u\,\right> \,=\, 0\; \;\forall\; u \,\in\, V\,\right\}$.
\item[$(ii)$]\,$\left<\,V\,\right>$\, be the right \,$\mathfrak{Q}$-linear subspace of \,$\mathbb{H}^{\,R}\,(\,\mathfrak{Q}\,)$\, consisting of all finite right \,$\mathfrak{Q}$-linear combinations of elements of \,$V$.  
\end{itemize}
\end{definition}

\begin{definition}\cite{RGV}
Every quaternionic Hilbert space \,$\mathbb{H}^{\,R}\,(\,\mathfrak{Q}\,)$\, admits a subset \,$N$, called Hilbert basis or orthonormal basis of \,$\mathbb{H}^{\,R}\,(\,\mathfrak{Q}\,)$\,, such that for \,$u,\, v \,\in\, N$, \,$\left<\,u,\, v\,\right> \,=\, 0$\, if \,$u \,\neq\, v$\, and \,$\left<\,u,\, u\,\right> \,=\, 1$.  
\end{definition}

\begin{theorem}\cite{RGV}
Let \,$\mathbb{H}^{\,R}\,(\,\mathfrak{Q}\,)$\, be a right quaternionic Hilbert space and \,$N$\, be Hilbert basis of \,$\mathbb{H}^{\,R}\,(\,\mathfrak{Q}\,)$.\,Then the following conditions are equivalent:
\begin{itemize}
\item[$(i)$]For every \,$u,\, v \,\in\, \mathbb{H}^{\,R}\,(\,\mathfrak{Q}\,)$, the series \,$\sum\limits_{z \,\in\, N}\,\left<\,u,\, z\,\right>\,\left<\,z,\, v\,\right>$\, converges absolutely and \,$\left<\,u,\, v\,\right> \,=\, \sum\limits_{z \,\in\, N}\,\left<\,u,\, z\,\right>\,\left<\,z,\, v\,\right>$.
\item[$(ii)$]For every \,$u \,\in\, \mathbb{H}^{\,R}\,(\,\mathfrak{Q}\,)$, \,$\|\,u\,\|^{\,2} \,=\, \sum\limits_{z \,\in\, N}\,\left|\,\left<\,z,\, u\,\right>\,\right|^{\,2}$.
\item[$(iii)$]\,\,$V^{\,\bot} \,=\, 0$.
\item[$(iii)$]\,\,$\left<\,N\,\right>$\, is dense in \,$\mathbb{H}^{\,R}\,(\,\mathfrak{Q}\,)$. 
\end{itemize} 
\end{theorem}

\begin{definition}\cite{SLA}
Let \,$\mathbb{H}^{\,R}\,(\,\mathfrak{Q}\,)$\, be a right quaternionic Hilbert space and \,$T$\, be an operator on \,$\mathbb{H}^{\,R}\,(\,\mathfrak{Q}\,)$.\,Then \,$T$\, is said to be right $\mathfrak{Q}$-linear if 
\[T\,\left(\,u\,\alpha \,+\, v\,\beta\,\right) \,=\, \alpha\,T\,(\,u\,) \,+\, \beta\,T\,(\,v\,)\,,\;\; \,\text{for all}\; u,\, v \,\in\, \mathbb{H}^{\,R}\,(\,\mathfrak{Q}\,)\; \;\text{and}\; \alpha,\, \beta \,\in\, \mathfrak{Q}.\]
\,$T$\, is said to be bounded if there exist \,$K \,>\, 0$\, such that \,$\left\|\,T\,(\,v\,)\,\right\| \,\leq\, K\,\|\,v\,\|$, for all \,$v \,\in\, \mathbb{H}^{\,R}\,(\,\mathfrak{Q}\,)$.\,The adjoint operator \,$T^{\,\ast}$\, of \,$T$\, is defined as \,$\left<\,v,\, T\,u\,\right> \,=\, \left<\,T^{\,\ast}\,v,\, u\,\right>$, for all \,$u,\, v \,\in\, \mathbb{H}^{\,R}\,(\,\mathfrak{Q}\,)$\, and \,$T$\, is said to be self-adjoint if \,$T \,=\, T^{\,\ast}$    
\end{definition}

\begin{theorem}\cite{SLA}
Let \,$\mathbb{H}^{\,R}\,(\,\mathfrak{Q}\,)$\, be a right quaternionic Hilbert space and \,$S,\, T$\, be two bounded right linear operators on \,$\mathbb{H}^{\,R}\,(\,\mathfrak{Q}\,)$.\,Then
\begin{itemize}
\item[$(i)$]\,$T \,+\, S$\, and \,$T\,S$\, are bounded right linear operators on \,$\mathbb{H}^{\,R}\,(\,\mathfrak{Q}\,)$.\,Furthermore \,$\|\,T \,+\, S\,\| \,\leq\, \|\,T\,\| \,+\, \|\,S\,\|$\, and \,$\|\,T\,S\,\| \,\leq\, \|\,T\,\|\,\|\,S\,\|$.
\item[$(ii)$]$(\,T \,+\, S\,)^{\,\ast} \,=\, T^{\,\ast} \,+\, S^{\,\ast}$, \,$(\,T\,S\,)^{\,\ast} \,=\, S^{\,\ast}\,T^{\,\ast}$\, and \,$(\,T^{\,\ast}\,)^{\,\ast} \,=\, T$.
\item[$(iii)$]\,$I_{H}^{\,\ast} \,=\, I_{H}$, where \,$I_{H}$\, is an identity operator on \,$\mathbb{H}^{\,R}\,(\,\mathfrak{Q}\,)$.
\item[$(iv)$]If \,$T$\, ia an invertible operator then \,$\left(\,T^{\,-\, 1}\,\right)^{\,\ast} \,=\, \left(\,T^{\,\ast}\,\right)^{\,-\, 1}$.  
\end{itemize}
\end{theorem}

\begin{theorem}\cite{RGV}
Let \,$\mathbb{H}^{\,R}\,(\,\mathfrak{Q}\,)$\, be a right quaternionic Hilbert space and let \,$T \,\in\, \mathcal{B}\left(\,\mathbb{H}^{\,R}\,(\,\mathfrak{Q}\,)\,\right)$\, be an operator.\,If \,$T \,\geq\, 0$, then exists a unique operator in \,$\mathcal{B}\left(\,\mathbb{H}^{\,R}\,(\,\mathfrak{Q}\,)\,\right)$, say \,$\sqrt{T}$, such that \,$\sqrt{T} \,\geq\, 0$\, and \,$\sqrt{T}\,\sqrt{T} \,=\, T$.\,Furthermore, \,$\sqrt{T}$\, commutes with every operator which commutes with \,$T$\, and if \,$T$\, is invertible and self-adjoint, then \,$\sqrt{T}$ is also invertible and self-adjoint. 
\end{theorem}   

Throughout this paper, \,$\mathfrak{Q}$\, is considered to be a non-commutative field of quaternions, \,$J$\, is subset of integers \,$\mathbb{Z}$\, and \,$\mathbb{H}^{\,R}\,(\,\mathfrak{Q}\,)$\, is a separable right quaternionic Hilbert space.\,By the term "right linear operator" we mean a "right \,$\mathfrak{Q}$-linear operator"   and \,$\mathcal{B}\left(\,\mathbb{H}^{\,R}\,(\,\mathfrak{Q}\,)\,\right)$\, denotes the set of all bounded (\,right \,$\mathfrak{Q}$-linear\,) operators on \,$\mathbb{H}^{\,R}\,(\,\mathfrak{Q}\,)$.

\section{Various generalizations of frame in quaternionic Hilbert space}

\begin{definition}\cite{SKS}
Let \,$\mathbb{H}^{\,R}\,(\,\mathfrak{Q}\,)$\, be a right quaternionic Hilbert space and \,$\left\{\,f_{\,j}\,\right\}_{j \,\in\, J}$\, be a sequence in \,$\mathbb{H}^{\,R}\,(\,\mathfrak{Q}\,)$.\,Then \,$\left\{\,f_{\,j}\,\right\}_{j \,\in\, J}$\, is a frame for \,$\mathbb{H}^{\,R}\,(\,\mathfrak{Q}\,)$\, if there exist constants \,$A,\, B \,>\, 0$\, such that
\[ A\, \|\,f\,\|^{\,2} \,\leq\, \sum\limits_{j \,\in\, J}\, \left|\, \left<\,f,\, f_{\,j} \, \right>\,\right|^{\,2} \,\leq\, B \,\|\,f\,\|^{\,2}\; \;\forall\; f \,\in\, \,\mathbb{H}^{\,R}\,(\,\mathfrak{Q}\,). \]
The constants \,$A$\, and \,$B$\, are called frame bounds.
\end{definition}

\begin{example}
Let \,$N$\, be a Hilbert basis for right separable quaternionic Hilbert space \,$\mathbb{H}^{\,R}\,(\,\mathfrak{Q}\,)$\, such that for \,$z_{\,i},\, z_{\,k} \,\in\, N$, \,$i,\, k \,\in\, J$, we have \,$\left<\,z_{\,i},\, z_{\,k}\,\right> \,=\, 0$\, if \,$i \,\neq\, k$\, and \,$\left<\,z_{\,i},\, z_{\,i}\,\right> \,=\, 1$.\,Let \,$\left\{\,f_{\,j}\,\right\}_{j \,\in\, J}$\, be a sequence in \,$\mathbb{H}^{\,R}\,(\,\mathfrak{Q}\,)$\, such that \,$u_{\,j} \,=\, u_{\,j \,+\, 1} \,=\, z_{\,j}\,,\; \, j \,\in\, J$.\,Then  \,$\left\{\,f_{\,j}\,\right\}_{j \,\in\, J}$\, is a tight frame for \,$\mathbb{H}^{\,R}\,(\,\mathfrak{Q}\,)$\, with bound \,$2$.    
\end{example}

\begin{definition}\cite{AMJ}
Let \,$\left\{\,W^{\,R}_{j}\,\right\}_{ j \,\in\, J}$\; be a collection of closed subspaces of a right separable quaternionic Hilbert space \,$\mathbb{H}^{\,R}\,(\,\mathfrak{Q}\,)$\, and \,$\left\{\,v_{j}\,\right\}_{ j \,\in\, J}$\, be a collection of positive weights.\;A family of weighted closed subspaces \,$\left\{\, (\,W^{\,R}_{j},\, v_{j}\,) \,:\, j \,\in\, J\,\right\}$\; is called a fusion frame for \,$\mathbb{H}^{\,R}\,(\,\mathfrak{Q}\,)$\, if there exist constants \,$0 \,<\, A \,\leq\, B \,<\, \infty$\; such that
\[A \;\left\|\,f \,\right\|^{\,2} \,\leq\, \sum\limits_{\,j \,\in\, J}\, v_{j}^{\,2}\,\left\|\,P_{\,W^{\,R}_{j}}\,(\,f\,) \,\right\|^{\,2} \,\leq\, B \; \left\|\,f \, \right\|^{\,2}\; \;\forall\; f \,\in\, \mathbb{H}^{\,R}\,(\,\mathfrak{Q}\,).\]
The constants \,$A,\, B$\; are called fusion frame bounds.\;If \,$A \,=\, B$\; then the fusion frame is called a tight fusion frame, if \,$A \,=\, B \,=\, 1$\; then it is called a Parseval fusion frame.
\end{definition}

\begin{definition}\cite{AMJ}
Let \,$\mathbb{H}^{\,R}\,(\,\mathfrak{Q}\,)$\, be a right quaternionic Hilbert space and \,$\left\{\,\mathbb{H}_{\,j}^{\,R}\,(\,\mathfrak{Q}\,)\,\right\}_{j \,\in\, J}$\, be a collection of right quaternionic Hilbert spaces.\,Then the sequence of bounded right linear operator \,$\left\{\,\Lambda_{j} \,\in\, \mathcal{B}\,(\,\mathbb{H}^{\,R}\,(\,\mathfrak{Q}\,),\, \mathbb{H}_{\,j}^{\,R}\,(\,\mathfrak{Q}\,)\,) \,:\, j \,\in\, J\,\right\}$\, is called frame of operator for \,$\mathbb{H}^{\,R}\,(\,\mathfrak{Q}\,)$\; with respect to \,$\left\{\,\mathbb{H}_{\,j}^{\,R}\,(\,\mathfrak{Q}\,)\,\right\}_{j \,\in\, J}$\, if there are two positive constants \,$A$\, and \,$B$\, such that
\[A \;\left \|\, f \,\right \|^{\,2} \,\leq\, \sum\limits_{\,j \,\in\, J}\, \left\|\,\Lambda_{j}\,f \,\right\|^{\,2} \,\leq\, B \; \left\|\, f \, \right\|^{\,2}\; \;\forall\; f \,\in\, \mathbb{H}^{\,R}\,(\,\mathfrak{Q}\,).\]
The constants \,$A$\; and \,$B$\; are called the lower and upper frame bounds, respectively.
\end{definition}

\begin{example}
Let \,$\mathbb{H}^{\,R}\,(\,\mathfrak{Q}\,)$\, be a right quaternionic Hilbert space and \,$N \,=\, \left\{\,u_{\,j}\,\right\}_{j \,\in\, \mathbb{N}}$\, be an orthonormal basis for \,$\mathbb{H}^{\,R}\,(\,\mathfrak{Q}\,)$.\,Define \,$\Lambda_{j} \,:\, \mathbb{H}^{\,R}\,(\,\mathfrak{Q}\,) \,\to\, \mathfrak{Q}$\, by 
\,$\Lambda_{j}\,(\,f\,) \,=\, \left<\,z_{\,j},\, f\,\right>$, for all \,$f \,\in\, \mathbb{H}^{\,R}\,(\,\mathfrak{Q}\,)$, \,$j \,\in\, \mathbb{N}$.\,Then \,$\left\{\,\Lambda_{j} \,:\, \mathbb{H}^{\,R}\,(\,\mathfrak{Q}\,) \,\to\, \mathfrak{Q}\,\right\}_{j \,\in\, \mathbb{N}}$\, is frame of operator for \,$\mathbb{H}^{\,R}\,(\,\mathfrak{Q}\,)$\, with respect to \,$\mathfrak{Q}$.  
\end{example}

\section{$g$-fusion frame in Quaternionic Hilbert space}

\smallskip\hspace{.6 cm} In this section, we present the concept of generalized fusion frame or \,$g$-fusion frame in a right quaternionic Hilbert space and discuss some few properties.

\begin{definition}
Let \,$W \,=\, \left\{\,W^{\,R}_{j}\,\right\}_{ j \,\in\, J}$\, be a collection of closed subspaces of right quaternionic Hilbert space \,$\mathbb{H}^{\,R}\,(\,\mathfrak{Q}\,)$\, and \,$\left\{\,v_{j}\,\right\}_{ j \,\in\, J}$\; be a collection of positive weights and \,$\left\{\,\Lambda_{j} \,:\, \mathbb{H}^{\,R}\,(\,\mathfrak{Q}\,) \,\to\, \mathbb{H}_{\,j}^{\,R}\,(\,\mathfrak{Q}\,)\,\right\}$\; be a collection of bounded right linear operators.\,Then the family \,$\Lambda \,=\, \left\{\,\left(\,W^{\,R}_{j},\, \Lambda_{j},\, v_{j}\,\right)\,\right\}_{j \,\in\, J}$\, is called a generalized fusion frame or a $g$-fusion frame for \,$\mathbb{H}^{\,R}\,(\,\mathfrak{Q}\,)$\, with respect to \,$\left\{\,\mathbb{H}_{\,j}^{\,R}\,(\,\mathfrak{Q}\,)\,\right\}_{j \,\in\, J}$\, if there exist constants \,$0 \,<\, A \,\leq\, B \,<\, \infty$\, such that
\begin{equation}\label{1.eq1.1}
A \;\left \|\,f \,\right \|^{\,2} \,\leq\, \sum\limits_{\,j \,\in\, J}\,v_{j}^{\,2}\, \left\|\,\Lambda_{j}\,P_{\,W^{\,R}_{j}}\,(\,f\,) \,\right\|^{\,2} \,\leq\, B \; \left\|\, f \, \right\|^{\,2}\; \;\forall\; f \,\in\, \mathbb{H}^{\,R}\,(\,\mathfrak{Q}\,),
\end{equation}
where each \,$P_{\,W^{\,R}_{j}}$\, is an orthogonal projection onto the closed subspace \,$W^{\,R}_{j}$\, and \,$\left\{\,\mathbb{H}_{\,j}^{\,R}\,(\,\mathfrak{Q}\,)\,\right\}_{j \,\in\, J}$\, is the collection of right quaternionic Hilbert spaces.\,The constants \,$A$\, and \,$B$\, are called the lower and upper bounds of $g$-fusion frame, respectively.\,If \,$A \,=\, B$\, then \,$\Lambda$\, is called tight $g$-fusion frame for \,$\mathbb{H}^{\,R}\,(\,\mathfrak{Q}\,)$\, and if \,$A \,=\, B \,=\, 1$\, then we say \,$\Lambda$\, is a Parseval $g$-fusion frame for \,$\mathbb{H}^{\,R}\,(\,\mathfrak{Q}\,)$.\,If  \,$\Lambda$\; satisfies only the right inequality of (\ref{1.eq1.1}) then it is called a $g$-fusion Bessel sequence with bound \,$B$\, in \,$\mathbb{H}^{\,R}\,(\,\mathfrak{Q}\,)$. 
\end{definition}

\begin{example}
Let \,$\mathbb{H}^{\,R}\,(\,\mathfrak{Q}\,)$\, be a right separable quaternionic Hilbert space and \,$\left\{\,z_{\,j}\,\right\}_{j \,\in\, \mathbb{N}}$\, be an orthonormal basis for \,$\mathbb{H}^{\,R}\,(\,\mathfrak{Q}\,)$.\,Define \[W^{\,R}_{1} \,=\, \mathbb{H}_{\,1}^{\,R}\,(\,\mathfrak{Q}\,) \,=\, \text{span}\left\{\,z_{\,1}\,\right\},\, \;W^{\,R}_{j} \,=\, \mathbb{H}_{\,j}^{\,R}\,(\,\mathfrak{Q}\,) \,=\, \text{span}\left\{\,z_{\,j \,-\, 1}\,\right\}\,,\; \,j \,\geq\, 2\]
and \,$v_{\,j} \,=\, 1$, \,$j \,\in\, \mathbb{N}$.\,Now, for each \,$j \,\in\, \mathbb{N}$, define \,$\Lambda_{j} \,:\, \mathbb{H}^{\,R}\,(\,\mathfrak{Q}\,) \,\to\, \mathbb{H}_{\,j}^{\,R}\,(\,\mathfrak{Q}\,)$\, by
\[\Lambda_{j}\,f \,=\, \left<\,z_{\,j},\, f\,\right>\,z_{\,j}\,,\; \,\text{for all}\; \;f \,\in\, \mathbb{H}^{\,R}\,(\,\mathfrak{Q}\,). \]
Then it is easy to verify that
\[\left \|\,f \,\right \|^{\,2} \,\leq\, \sum\limits_{\,j \,\in\, \mathbb{N}}\,\left\|\,\Lambda_{j}\,P_{\,W^{\,R}_{j}}\,(\,f\,) \,\right\|^{\,2} \,\leq\, 2\,\left\|\, f \, \right\|^{\,2}\; \;\forall\; f \,\in\, \mathbb{H}^{\,R}\,(\,\mathfrak{Q}\,).\]
Thus, \,$\Lambda \,=\, \left\{\,\left(\,W^{\,R}_{j},\, \Lambda_{j},\, 1\,\right)\,\right\}_{j \,\in\, \mathbb{N}}$\, is a \,$g$-fusion frame for \,$\mathbb{H}^{\,R}\,(\,\mathfrak{Q}\,)$\, with bounds \,$1$\, and \,$2$.    
\end{example}

Define the space
\[\mathcal{H}_{2}  \,=\, \bigoplus_{j \,\in\, J}\,\mathbb{H}_{\,j}^{\,R}\,(\,\mathfrak{Q}\,) \,=\,  \left \{\,\{\,f_{\,j}\,\}_{j \,\in\, J} \,:\, f_{\,j} \,\in\, \mathbb{H}_{\,j}^{\,R}\,(\,\mathfrak{Q}\,),\; \sum\limits_{\,j \,\in\, J}\, \left \|\,f_{\,j}\,\right \|_{\,\mathbb{H}_{\,j}^{\,R}\,(\,\mathfrak{Q}\,)}^{\,2} \,<\, \infty \,\right\}\]
under right multiplications by quaternionic scalars together with the quaternionic inner product is given by \[\left<\,\{\,f_{\,j}\,\}_{ j \,\in\, J},\, \{\,g_{\,j}\,\}_{ j \,\in\, J}\,\right> \;=\; \sum\limits_{\,j \,\in\, J}\, \left<\,f_{\,j},\, g_{\,j}\,\right>_{\mathbb{H}_{\,j}^{\,R}\,(\,\mathfrak{Q}\,)},\]
and the norm is defined as \,$\left\|\,\{\,f_{\,j}\,\}_{ j \,\in\, J}\,\right\|_{\mathcal{H}_{2}} \,=\, \sum\limits_{\,j \,\in\, J}\,\left\|\,f_{\,j}\,\right\|_{\mathbb{H}_{\,j}^{\,R}\,(\,\mathfrak{Q}\,)}$, for all \,$\{\,f_{\,j}\,\}_{j \,\in\, J} \,\in\, \mathcal{H}_{2}$.\,It is easy to very that \,$\mathcal{H}_{2}$\, is a right quaternionic Hilbert space with respect to the quaternionic inner product given by above.

\begin{note}
Let \,$\Lambda$\, be \,$g$-fusion Bessel sequence for \,$\mathbb{H}^{\,R}\,(\,\mathfrak{Q}\,)$\, with bound \,$B$.\,Then for every sequence \,$\{\,f_{\,j}\,\}_{j \,\in\, J} \,\in\, \mathcal{H}_{2}$, the series \,$\sum\limits_{\,j \,\in\, J}\,v_{j}\,P_{\,W^{\,R}_{j}}\,\Lambda^{\,\ast}_{j}\,f_{\,j}$\, converges unconditionally.  
\end{note}

\begin{theorem}\label{1.th1.10}
The family \,$\Lambda$\, is a \,$g$-fusion Bessel sequence for \,$\mathbb{H}^{\,R}\,(\,\mathfrak{Q}\,)$\, with bound \,$B$\, if and only if the right linear operator \,$T_{\mathfrak{Q}} \,:\, \mathcal{H}_{2} \,\to\, \mathbb{H}^{\,R}\,(\,\mathfrak{Q}\,)$\, defined by
\[T_{\mathfrak{Q}}\,\left(\,\left\{\,f_{\,j}\,\right\}_{j \,\in\, J}\,\right) \,=\,  \sum\limits_{\,j \,\in\, J}\, v_{j}\, P_{\,W^{\,R}_{j}}\,\Lambda_{j}^{\,\ast}\,f_{j}\; \;\;\forall\; \{\,f_{j}\,\}_{j \,\in\, J} \,\in\, \mathcal{H}_{2},\]
is a well-defined and bounded operator with \,$\left\|\,T_{\mathfrak{Q}}\,\right\| \,\leq\, \sqrt{\,B}$.  
\end{theorem}

\begin{proof}
Suppose \,$\Lambda$\, is a \,$g$-fusion Bessel sequence for \,$\mathbb{H}^{\,R}\,(\,\mathfrak{Q}\,)$\, with bound \,$B$.\,Let \,$I$\, be a finite subset of \,$J$.\,Then
\begin{align*}
\left\|\,\sum\limits_{\,j \,\in\, I}\, v_{j}\, P_{\,W^{\,R}_{j}}\,\Lambda_{j}^{\,\ast}\,f_{j}\,\right\|^{\,2}& \,=\, \sup\limits_{\|\,g\,\| \,=\, 1}\,\left|\,\left<\,\sum\limits_{\,j \,\in\, I}\, v_{j}\, P_{\,W^{\,R}_{j}}\,\Lambda_{j}^{\,\ast}\,f_{j}\,,\, g\,\right>\,\right|^{\,2} \\
&=\,\sup\limits_{\|\,g\,\| \,=\, 1}\,\sum\limits_{\,j \,\in\, I}\,\left|\,\left<\,f_{j}\,,\, v_{j}\,\Lambda_{j}\,P_{\,W^{\,R}_{j}}\,(\,g\,)\,\right>\,\right|^{\,2}\\
&\leq\,\sum\limits_{\,j \,\in\, I}\,\left\|\,f_{\,j}\,\right\|^{\,2}_{\mathbb{H}_{\,j}^{\,R}\,(\,\mathfrak{Q}\,)}\,\sup\limits_{\|\,g\,\| \,=\, 1}\,\sum\limits_{\,j \,\in\, I}\,v_{j}^{\,2}\, \left\|\,\Lambda_{j}\,P_{\,W^{\,R}_{j}}\,(\,g\,) \,\right\|^{\,2}\\
&\leq\, B\,\sum\limits_{\,j \,\in\, I}\,\left\|\,f_{\,j}\,\right\|^{\,2}_{\mathbb{H}_{\,j}^{\,R}\,(\,\mathfrak{Q}\,)} \,<\, \infty.
\end{align*}
Thus, the series \,$\sum\limits_{\,j \,\in\, J}\,v_{j}\,P_{\,W^{\,R}_{j}}\,\Lambda^{\,\ast}_{j}\,f_{\,j}$\, converges unconditionally.\,Hence, the right linear operator \,$T_{\mathfrak{Q}}$\, is well-defined.\,By the above similar calculation it is easy to verify that \,$T_{\mathfrak{Q}}$\, is bounded and \,$\left\|\,T_{\mathfrak{Q}}\,\right\| \,\leq\, \sqrt{\,B}$. \\ 

Conversely, suppose that \,$T_{\mathfrak{Q}}$\, is well-defined and bounded right linear operator with \,$\left\|\,T_{\mathfrak{Q}}\,\right\| \,\leq\, \sqrt{\,B}$.\,Then the adjoint \,$T^{\,\ast}_{\mathfrak{Q}}$\, of a bounded right linear operator \,$T_{\mathfrak{Q}}$\, is itself bounded and \,$\left\|\,T_{\mathfrak{Q}}\,\right\| \,=\, \left\|\,T^{\,\ast}_{\mathfrak{Q}}\,\right\|$.\,Now, for \,$f \,\in\, \mathbb{H}^{\,R}\,(\,\mathfrak{Q}\,)$, we have
\begin{align*}
&\sum\limits_{\,j \,\in\, J}\,v_{j}^{\,2}\, \left\|\,\Lambda_{j}\,P_{\,W^{\,R}_{j}}\,(\,f\,) \,\right\|^{\,2} \,=\, \left\|\,T^{\,\ast}_{\mathfrak{Q}}\,f\,\right\|^{\,2} \,\leq\, \left\|\,T_{\mathfrak{Q}}\,\right\|^{\,2}\,\|\,f\,\|^{\,2}  \,\leq\, B\,\|\,f\,\|^{\,2}.
\end{align*} 
Thus, \,$\Lambda$\, is a \,$g$-fusion Bessel sequence for the right quaternionic Hilbert space \,$\mathbb{H}^{\,R}\,(\,\mathfrak{Q}\,)$\, with bound \,$B$.          
\end{proof} 

Let \,$\Lambda$\, be a \,$g$-fusion Bessel sequence for the right quaternionic Hilbert space \,$\mathbb{H}^{\,R}\,(\,\mathfrak{Q}\,)$.\,Then the right linear operator \,$T_{\mathfrak{Q}} \,:\, \mathcal{H}_{2} \,\to\, \mathbb{H}^{\,R}\,(\,\mathfrak{Q}\,)$\, given by
\[T_{\mathfrak{Q}}\,\left(\,\left\{\,f_{\,j}\,\right\}_{j \,\in\, J}\,\right) \,=\,  \sum\limits_{\,j \,\in\, J}\, v_{j}\, P_{\,W^{\,R}_{j}}\,\Lambda_{j}^{\,\ast}\,f_{j}\; \;\;\forall\; \{\,f_{j}\,\}_{j \,\in\, J} \,\in\, \mathcal{H}_{2},\]
is called the (\,right\,) synthesis operator and the adjoint of \,$T_{\mathfrak{Q}}$\, given by 
\[T^{\,\ast}_{\mathfrak{Q}} \,:\, \mathbb{H}^{\,R}\,(\,\mathfrak{Q}\,) \,\to\, \mathcal{H}_{2}\;,\; \,T^{\,\ast}_{\mathfrak{Q}}\,(\,f\,) \,=\,  \left\{\,v_{j}\,\Lambda_{j}\, P_{\,W^{\,R}_{j}}\,(\,f\,)\,\right\}_{ j \,\in\, J}\; \;\forall\; f \,\in\, H,\]
is called the (\,right\,) analysis operator.

\begin{definition}
Let \,$\Lambda$\, be a \,$g$-fusion frame for the right quaternionic Hilbert space \,$\mathbb{H}^{\,R}\,(\,\mathfrak{Q}\,)$.\,The right linear operator \,$S_{\mathfrak{Q}} \,:\, \mathbb{H}^{\,R}\,(\,\mathfrak{Q}\,) \,\to\, \mathbb{H}^{\,R}\,(\,\mathfrak{Q}\,)$\, defined by
\[S_{\mathfrak{Q}}\,f \,=\, T_{\mathfrak{Q}}\,T^{\,\ast}_{\mathfrak{Q}}\,f \,=\, \sum\limits_{\,j \,\in\, J}\, v^{\,2}_{j}\, P_{\,W^{\,R}_{j}}\,\Lambda_{j}^{\,\ast}\,\Lambda_{j}\,P_{\,W^{\,R}_{j}}\,f\;,\;\, f \,\in\, \mathbb{H}^{\,R}\,(\,\mathfrak{Q}\,),\] 
is called the (\,right\,) \,$g$-fusion frame operator for \,$\Lambda$.  
\end{definition}

In the next Theorem, we will discuss a few properties of the frame operator for the \,$g$-fusion frame in right quaternionic Hilbert space.

\begin{theorem}
Let \,$\Lambda$\, be a \,$g$-fusion frame for \,$\mathbb{H}^{\,R}\,(\,\mathfrak{Q}\,)$\, with bounds \,$A,\,B$\, and \,$S_{\mathfrak{Q}}$\, be the corresponding right \,$g$-fusion frame operator.\,Then \,$S_{\mathfrak{Q}}$\, is positive, bounded, invertible and self-adjoint right linear operator on \,$\mathbb{H}^{\,R}\,(\,\mathfrak{Q}\,)$.  
\end{theorem} 

\begin{proof}
For each \,$f \,\in\, \mathbb{H}^{\,R}\,(\,\mathfrak{Q}\,)$, we have \,$\left<\,S_{\mathfrak{Q}}\,f,\, f\,\right>  \,=\, \sum\limits_{\,j \,\in\, J}\,v_{j}^{\,2}\, \left\|\,\Lambda_{j}\,P_{\,W^{\,R}_{j}}\,(\,f\,) \,\right\|^{\,2}$\, and from (\ref{1.eq1.1}), we get
\begin{align*}
A \;\left \|\,f \,\right \|^{\,2} \,\leq\, \left<\,S_{\mathfrak{Q}}\,f,\, f\,\right>  \,\leq\, B \; \left\|\, f \, \right\|^{\,2} \,\Rightarrow\, A\,I_{H} \,\leq\, S_{\mathfrak{Q}} \,\leq\, B\,I_{H}. 
\end{align*}
Hence, \,$S_{\mathfrak{Q}}$\, is positive and bounded right linear operator on \,$\mathbb{H}^{\,R}\,(\,\mathfrak{Q}\,)$\, and consequently it is a invertible.

Furthermore, for any \,$f,\, g \,\in\, \mathbb{H}^{\,R}\,(\,\mathfrak{Q}\,)$, we have
\begin{align*}
\left<\,S_{\mathfrak{Q}}\,f,\, g\,\right> &\,=\, \left<\,\sum\limits_{\,j \,\in\, J}\, v^{\,2}_{j}\, P_{\,W^{\,R}_{j}}\,\Lambda_{j}^{\,\ast}\,\Lambda_{j}\,P_{\,W^{\,R}_{j}}\,f,\, g\,\right> \\
&\,=\, \sum\limits_{\,j \,\in\, J}\,\left<\,f,\, v^{\,2}_{j}\, P_{\,W^{\,R}_{j}}\,\Lambda_{j}^{\,\ast}\,\Lambda_{j}\,P_{\,W^{\,R}_{j}}\,g\,\right> \\
&=\,\left<\,f,\, \sum\limits_{\,j \,\in\, J}\,v^{\,2}_{j}\,P_{\,W^{\,R}_{j}}\,\Lambda_{j}^{\,\ast}\,\Lambda_{j}\,P_{\,W^{\,R}_{j}}\,(\,g\,)\,\right> \,=\, \left<\,f,\, S_{\mathfrak{Q}}\,g\,\right>. 
\end{align*}
Thus, \,$S_{\mathfrak{Q}}$\, is also self-adjoint right linear operator on \,$\mathbb{H}^{\,R}\,(\,\mathfrak{Q}\,)$.       
\end{proof}

\begin{corollary}
For every \,$f \,\in\, \mathbb{H}^{\,R}\,(\,\mathfrak{Q}\,)$, we get the reconstruction formula as:
\[f \,=\, \sum\limits_{\,j \,\in\, J}\, v^{\,2}_{j}\,S^{\,-\, 1}_{\mathfrak{Q}}\,P_{\,W^{\,R}_{j}}\,\Lambda_{j}^{\,\ast}\,\Lambda_{j}\,P_{\,W^{\,R}_{j}}\,f \,=\, \sum\limits_{\,j \,\in\, J}\, v^{\,2}_{j}\,P_{\,W^{\,R}_{j}}\,\Lambda_{j}^{\,\ast}\,\Lambda_{j}\,P_{\,W^{\,R}_{j}}\,S^{\,-\, 1}_{\mathfrak{Q}}\,f.\]
\end{corollary}

In the following Theorem, we establish a characterization of a Parseval \,$g$-fusion frame for the right quaternionic Hilbert space \,$\mathbb{H}^{\,R}\,(\,\mathfrak{Q}\,)$.  

\begin{theorem}
Let \,$\Lambda$\, be a \,$g$-fusion frame for \,$\mathbb{H}^{\,R}\,(\,\mathfrak{Q}\,)$\, with the corresponding right \,$g$-fusion frame operator \,$S_{\mathfrak{Q}}$.\,Then \,$\Lambda$\, is a Parseval \,$g$-fusion frame for \,$\mathbb{H}^{\,R}\,(\,\mathfrak{Q}\,)$\, if and only if \,$S_{\mathfrak{Q}}$\, is an identity operator on \,$\mathbb{H}^{\,R}\,(\,\mathfrak{Q}\,)$.     
\end{theorem}

\begin{proof}
Let \,$\Lambda$\, be a Parseval \,$g$-fusion frame for \,$\mathbb{H}^{\,R}\,(\,\mathfrak{Q}\,)$.\,Then, for each \,$f \,\in\, \mathbb{H}^{\,R}\,(\,\mathfrak{Q}\,)$, we get
\begin{align*}
&\sum\limits_{\,j \,\in\, J}\,v_{j}^{\,2}\, \left\|\,\Lambda_{j}\,P_{\,W^{\,R}_{j}}\,(\,f\,) \,\right\|^{\,2} \,=\, \|\,f\,\|^{\,2} \,\Rightarrow\, \left<\,S_{\mathfrak{Q}}\,f,\, f\,\right> \,=\, \left<\,f,\, f\,\right>.
\end{align*}
This is shows that \,$S_{\mathfrak{Q}}$\, is an identity operator on \,$\mathbb{H}^{\,R}\,(\,\mathfrak{Q}\,)$.\\

Conversely, suppose that \,$S_{\mathfrak{Q}}$\, is an identity operator on \,$\mathbb{H}^{\,R}\,(\,\mathfrak{Q}\,)$.\,Then, for \,$f \,\in\, \mathbb{H}^{\,R}\,(\,\mathfrak{Q}\,)$, we get \,$f \,=\, S_{\mathfrak{Q}}\,f \,=\, \sum\limits_{\,j \,\in\, J}\, v^{\,2}_{j}\, P_{\,W^{\,R}_{j}}\,\Lambda_{j}^{\,\ast}\,\Lambda_{j}\,P_{\,W^{\,R}_{j}}\,f$.\,Therefore, for \,$f \,\in\, \mathbb{H}^{\,R}\,(\,\mathfrak{Q}\,)$, we have
\begin{align*}
\|\,f\,\|^{\,2}& \,=\, \left<\,f,\, f\,\right> \,=\, \left<\,\sum\limits_{\,j \,\in\, J}\, v^{\,2}_{j}\, P_{\,W^{\,R}_{j}}\,\Lambda_{j}^{\,\ast}\,\Lambda_{j}\,P_{\,W^{\,R}_{j}}\,f,\, f\,\right>\\
& =\,\sum\limits_{\,j \,\in\, J}\, v^{\,2}_{j}\,\left<\,\Lambda_{j}\,P_{\,W^{\,R}_{j}}\,f,\, \Lambda_{j}\,P_{\,W^{\,R}_{j}}\,f\,\right> \,=\, \sum\limits_{\,j \,\in\, J}\,v_{j}^{\,2}\, \left\|\,\Lambda_{j}\,P_{\,W^{\,R}_{j}}\,(\,f\,) \,\right\|^{\,2}. 
\end{align*}
Thus, \,$\Lambda$\, is a Parseval \,$g$-fusion frame for the right quaternionic Hilbert space \,$\mathbb{H}^{\,R}\,(\,\mathfrak{Q}\,)$.     
\end{proof}

Next, we give a characterization of a \,$g$-fusion frame for the right quaternionic Hilbert space \,$\mathbb{H}^{\,R}\,(\,\mathfrak{Q}\,)$\, with the help of its right synthesis operator.

\begin{theorem}
The family \,$\Lambda$\, is a \,$g$-fusion frame for \,$\mathbb{H}^{\,R}\,(\,\mathfrak{Q}\,)$\, if and only if the right synthesis operator \,$T_{\mathfrak{Q}}$\, is well-defined and bounded mapping from \,$\mathcal{H}_{\,2}$\, onto \,$\mathbb{H}^{\,R}\,(\,\mathfrak{Q}\,)$.  
\end{theorem}

\begin{proof}
Let \,$\Lambda$\, is a \,$g$-fusion frame for \,$\mathbb{H}^{\,R}\,(\,\mathfrak{Q}\,)$.\,Then it is easy to verify that \,$T_{\mathfrak{Q}}$\, is well-defined and bounded mapping from \,$\mathcal{H}_{\,2}$\, onto \,$\mathbb{H}^{\,R}\,(\,\mathfrak{Q}\,)$.\\

Conversely, suppose that the right synthesis operator \,$T_{\mathfrak{Q}}$\, is well-defined and bounded mapping from \,$\mathcal{H}_{\,2}$\, onto \,$\mathbb{H}^{\,R}\,(\,\mathfrak{Q}\,)$.\,Then by Theorem \ref{1.th1.10}, \,$\Lambda$\, is a \,$g$-fusion Bessel sequence for \,$\mathbb{H}^{\,R}\,(\,\mathfrak{Q}\,)$.\,Since \,$T_{\mathfrak{Q}}$\, is onto, there exists a right linear bounded operator \,$T^{\,\dagger}_{\mathfrak{Q}} \,:\, \mathbb{H}^{\,R}\,(\,\mathfrak{Q}\,) \,\to\, \mathcal{H}_{2}$\, such that
\[f \,=\, T_{\mathfrak{Q}}\,T^{\,\dagger}_{\mathfrak{Q}}\,f \,=\, \sum\limits_{\,j \,\in\, J}\, v_{j}\, P_{\,W^{\,R}_{j}}\,\Lambda_{j}^{\,\ast}\,\left(\,T^{\,\dagger}_{\mathfrak{Q}}\,f\,\right)_{j}\,,\; \; f \,\in\, \mathbb{H}^{\,R}\,(\,\mathfrak{Q}\,),\]
where \,$\left(\,T^{\,\dagger}_{\mathfrak{Q}}\,f\,\right)_{j}$\, denotes the \,$j$-th coordinate of \,$T^{\,\dagger}_{\mathfrak{Q}}\,f$.\,Now, for each \,$f \,\in\, \mathbb{H}^{\,R}\,(\,\mathfrak{Q}\,)$, we have
\begin{align*}
\|\,f\,\|^{\,4}& \,=\, \left|\,\left<\,f,\, f\,\right>\,\right|^{\,2} \,=\, \left|\,\left<\,\sum\limits_{\,j \,\in\, J}\, v_{j}\, P_{\,W^{\,R}_{j}}\,\Lambda_{j}^{\,\ast}\,\left(\,T^{\,\dagger}_{\mathfrak{Q}}\,f\,\right)_{j},\, f\,\right>\,\right|^{\,2}\\
&\leq\,\sum\limits_{\,j \,\in\, J}\,\left|\,\left(\,T^{\,\dagger}_{\mathfrak{Q}}\,f\,\right)_{j}\,\right|^{\,2}\,\sum\limits_{\,j \,\in\, J}\,v_{j}^{\,2}\, \left\|\,\Lambda_{j}\,P_{\,W^{\,R}_{j}}\,(\,f\,) \,\right\|^{\,2}\\
&\leq\,\left\|\,T^{\,\dagger}_{\mathfrak{Q}}\,\right\|^{\,2}\,\|\,f\,\|^{\,2}\,\sum\limits_{\,j \,\in\, J}\,v_{j}^{\,2}\, \left\|\,\Lambda_{j}\,P_{\,W^{\,R}_{j}}\,(\,f\,) \,\right\|^{\,2}.\\
&\Rightarrow\, \dfrac{1}{\left\|\,T^{\,\dagger}_{\mathfrak{Q}}\,\right\|^{\,2}}\,\|\,f\,\|^{\,2} \,\leq\, \sum\limits_{\,j \,\in\, J}\,v_{j}^{\,2}\, \left\|\,\Lambda_{j}\,P_{\,W^{\,R}_{j}}\,(\,f\,) \,\right\|^{\,2}.  
\end{align*}
Thus, \,$\Lambda$\, is a \,$g$-fusion frame for the right quaternionic Hilbert space \,$\mathbb{H}^{\,R}\,(\,\mathfrak{Q}\,)$.     
\end{proof} 

\begin{note}\label{nt1.01}
Let \,$V \,\subset\, \mathbb{H}^{\,R}\,(\,\mathfrak{Q}\,)$\; be a closed subspace and \,$T \,\in\, \mathcal{B}\left(\,\mathbb{H}^{\,R}\,(\,\mathfrak{Q}\,)\,\right)$. Then \,$P_{\,V}\, T^{\,\ast} \,=\, P_{\,V}\,T^{\,\ast}\, P_{\,\overline{T\,V}}$.
\end{note}

In the next Theorem, we will construct a \,$g$-fusion frame with the help of a given \,$g$-fusion frame in a right quaternionic Hilbert space.

\begin{theorem}\label{1.th1.11}
Let \,$\Lambda \,=\, \{\,\left(\,W^{\,R}_{j},\, \Lambda_{j},\, v_{j}\,\right)\,\}_{j \,\in\, J}$\, be \,$g$-fusion frame for \,$\mathbb{H}^{\,R}\,(\,\mathfrak{Q}\,)$\, with bounds \,$A,\,B$\, and \,$S_{\mathfrak{Q}}$\, be the corresponding right \,$g$-fusion frame operator.\,Then \,$\Gamma \,=\, \left\{\,\left(\,S^{\,-\, 1}_{\mathfrak{Q}}\,W^{\,R}_{j},\, \Lambda_{j}\,P_{\,W^{\,R}_{j}}\,S^{\,-\, 1}_{\mathfrak{Q}},\, v_{j}\,\right)\,\right\}_{j \,\in\, J}$\, is a \,$g$-fusion frame for \,$\mathbb{H}^{\,R}\,(\,\mathfrak{Q}\,)$\, with bounds \,$1 \,/\, B$\, and \,$1 \,/\, A$.    
\end{theorem}

\begin{proof}
For each \,$f \,\in\, \mathbb{H}^{\,R}\,(\,\mathfrak{Q}\,)$, we have
\begin{align*}
\sum\limits_{\,j \,\in\, J}\,v_{j}^{\,2}\, \left\|\,\Lambda_{j}\,P_{\,W^{\,R}_{j}}\,S^{\,-\, 1}_{\mathfrak{Q}}\,P_{\,S^{\,-\, 1}_{\mathfrak{Q}}\,W^{\,R}_{j}}\,(\,f\,) \,\right\|^{\,2}& \,=\, \sum\limits_{\,j \,\in\, J}\,v_{j}^{\,2}\,\left\|\,\Lambda_{j}\,P_{\,W^{\,R}_{j}}\,S^{\,-\, 1}_{\mathfrak{Q}}\,(\,f\,)\,\right\|^{\,2}\\
&\leq\,B\,\left\|\,S^{\,-\, 1}_{\mathfrak{Q}}\,\right\|^{\,2}\,\|\,f\,\|^{\,2}. 
\end{align*} 
Thus, \,$\Gamma$\, is a \,$g$-fusion Bessel sequence for \,$\mathbb{H}^{\,R}\,(\,\mathfrak{Q}\,)$.\,So, the right \,$g$-fusion frame operator for \,$\Gamma$\, is well-defined.\,Now, it is easy to verify that the right \,$g$-fusion frame operator for \,$\Gamma$\, is \,$S^{\,-\, 1}_{\mathfrak{Q}}$.\,The operator \,$S^{\,-\, 1}_{\mathfrak{Q}}$\, commutes with both \,$S_{\mathfrak{Q}}$\, and \,$I_{H}$.\,Thus, multiplying the inequality \,$A\,I_{H} \,\leq\, S_{\mathfrak{Q}} \,\leq\, B\,I_{H}$\, with \,$S^{\,-\, 1}_{\mathfrak{Q}}$, we get
\begin{align}
&B^{\,-\, 1}\,I_{H} \,\leq\, S^{\,-\, 1}_{\mathfrak{Q}} \,\leq\, A^{\,-\, 1}\,I_{H}\nonumber \\
&\Rightarrow\, B^{\,-\, 1} \;\left \|\,f \,\right \|^{\,2} \,\leq\, \left<\,S^{\,-\, 1}_{\mathfrak{Q}}\,f,\, f\,\right>  \,\leq\, A^{\,-\, 1} \; \left\|\, f \, \right\|^{\,2}\,,\; \,f \,\in\, \mathbb{H}^{\,R}\,(\,\mathfrak{Q}\,).\label{1.eq1.12}
\end{align} 
Now, \,$f \,\in\, \mathbb{H}^{\,R}\,(\,\mathfrak{Q}\,)$, we have
\begin{align*}
S^{\,-\, 1}_{\mathfrak{Q}}\,f &\,=\, S^{\,-\, 1}_{\mathfrak{Q}}\,\left(\,S_{\mathfrak{Q}}\,\left(\,S^{\,-\, 1}_{\mathfrak{Q}}\,f\,\right)\,\right)\\
&=\,S^{\,-\, 1}_{\mathfrak{Q}}\left(\,\sum\limits_{\,j \,\in\, J}\, v^{\,2}_{j}\,P_{\,W^{\,R}_{j}}\,\Lambda_{j}^{\,\ast}\,\Lambda_{j}\,P_{\,W^{\,R}_{j}}\,S^{\,-\, 1}_{\mathfrak{Q}}\,f\,\right)\\
&=\,\sum\limits_{\,j \,\in\, J}\, v^{\,2}_{j}\,S^{\,-\, 1}_{\mathfrak{Q}}\,P_{\,W^{\,R}_{j}}\,\Lambda_{j}^{\,\ast}\,\Lambda_{j}\,P_{\,W^{\,R}_{j}}\,S^{\,-\, 1}_{\mathfrak{Q}}\,f\\
&=\,\sum\limits_{\,j \,\in\, J}\, v^{\,2}_{j}\,\left(\,\Lambda_{j}\,P_{\,W^{\,R}_{j}}\,S^{\,-\, 1}_{\mathfrak{Q}}\,\right)^{\,\ast}\Lambda_{j}\,P_{\,W^{\,R}_{j}}\,S^{\,-\, 1}_{\mathfrak{Q}}\,f\\
&=\,\sum\limits_{\,j \,\in\, J}\, v^{\,2}_{j}\,P_{\,S^{\,-\, 1}_{\mathfrak{Q}}\,W^{\,R}_{j}}\,S^{\,-\, 1}_{\mathfrak{Q}}\,P_{\,W^{\,R}_{j}}\,\Lambda_{j}^{\,\ast}\,\Lambda_{j}\,P_{\,W^{\,R}_{j}}\,S^{\,-\, 1}_{\mathfrak{Q}}\,P_{\,S^{\,-\, 1}_{\mathfrak{Q}}\,W^{\,R}_{j}}\,f.  
\end{align*} 
Therefore, from (\ref{1.eq1.12}), for \,$f \,\in\, \mathbb{H}^{\,R}\,(\,\mathfrak{Q}\,)$, we get
\[B^{\,-\, 1} \;\left \|\,f \,\right \|^{\,2} \,\leq\,  \sum\limits_{\,j \,\in\, J}\,v_{j}^{\,2}\, \left\|\,\Lambda_{j}\,P_{\,W^{\,R}_{j}}\,S^{\,-\, 1}_{\mathfrak{Q}}\,P_{\,S^{\,-\, 1}_{\mathfrak{Q}}\,W^{\,R}_{j}}\,(\,f\,) \,\right\|^{\,2} \,\leq\, A^{\,-\, 1} \; \left\|\, f \, \right\|^{\,2}.\] 
This shows that \,$\Gamma$\, is a \,$g$-fusion frame for \,$\mathbb{H}^{\,R}\,(\,\mathfrak{Q}\,)$\, with bounds \,$1 \,/\, B$\, and \,$1 \,/\, A$.   
\end{proof}

In Theorem \ref{1.th1.11}, the \,$g$-fusion frame \,$\Gamma$\, for \,$\mathbb{H}^{\,R}\,(\,\mathfrak{Q}\,)$\, is called the canonical dual of the \,$g$-fusion frame of \,$\Lambda$.

\begin{note}
Let \,$\Lambda$\, be a \,$g$-fusion frame for \,$\mathbb{H}^{\,R}\,(\,\mathfrak{Q}\,)$\, with the corresponding right \,$g$-fusion frame operator \,$S_{\mathfrak{Q}}$.\,Then \,$\left\{\,\left(\,S^{\,-\, 1 \,/\, 2}_{\mathfrak{Q}}\,W^{\,R}_{j},\, \Lambda_{j}\,P_{\,W^{\,R}_{j}}\,S^{\,-\, 1 \,/\, 2}_{\mathfrak{Q}},\, v_{j}\,\right)\,\right\}_{j \,\in\, J}$\, is a Parseval \,$g$-fusion frame for \,$\mathbb{H}^{\,R}\,(\,\mathfrak{Q}\,)$. 
\end{note}

\begin{theorem}
Let \,$U \,\in\, \mathcal{B}\,\left(\,\mathbb{H}^{\,R}\,(\,\mathfrak{Q}\,\right)\,)$\, be an invertible bounded right $\mathfrak{Q}$-linear operator on \,$\mathbb{H}^{\,R}\,(\,\mathfrak{Q}\,)$\, and \,$\Lambda \,=\, \left\{\,\left(\,W_{j},\, \Lambda_{j},\, v_{j}\,\right)\,\right\}_{j \,\in\, J}$\, be a \,$g$-fusion frame for \,$\mathbb{H}^{\,R}\,(\,\mathfrak{Q}\,)$\, with bounds \,$A$\, and \,$B$.\,Then the family \,$\Gamma \,=\, \left\{\,\left(\,U\,W_{j},\, \Lambda_{j}\,P_{\,W_{j}}\,U^{\,\ast},\, v_{j}\,\right)\,\right\}_{j \,\in\, J}$\; is a \,$U\,U^{\,\ast}$-$g$-fusion frame for \,$\mathbb{H}^{\,R}\,(\,\mathfrak{Q}\,)$. 
\end{theorem}

\begin{proof}
Since \,$U$\, is an bounded right $\mathfrak{Q}$-linear operator on \,$\mathbb{H}^{\,R}\,(\,\mathfrak{Q}\,)$\, for any \,$j \,\in\, J$, \,$U\,W_{j}$\, is closed in \,$\mathbb{H}^{\,R}\,(\,\mathfrak{Q}\,)$.\;Now, for each \,$f \,\in\, \mathbb{H}^{\,R}\,(\,\mathfrak{Q}\,)$, using Note \ref{nt1.01}, we obtain
\begin{align*}
\sum\limits_{\,j \,\in\, J}\, v_{j}^{\,2}\, \left\|\,\Lambda_{j}\,P_{\,W^{\,R}_{j}}\,U^{\,\ast}\,P_{\,U\,W^{\,R}_{j}}\,(\,f\,)\,\right\|^{\,2}& \,=\, \sum\limits_{\,j \,\in\, J}\, v_{j}^{\,2}\, \left\|\,\Lambda_{j}\,P_{\,W^{\,R}_{j}}\,U^{\,\ast}\,(\,f\,)\,\right\|^{\,2}\\
&\leq\, B\; \left\|\,U^{\,\ast}\,f\,\right\|^{\,2} \,\leq\, B\, \|\,U\,\|^{\,2}\, \|\,f\,\|^{\,2}.
\end{align*}
On the other hand, for each \,$f \,\in\, \mathbb{H}^{\,R}\,(\,\mathfrak{Q}\,)$, we get
\begin{align*}
\dfrac{A}{\|\,U\,\|^{\,2}}\, \left\|\,\left(\,U\,U^{\,\ast}\,\right)^{\,\ast}\,f\,\right\|^{\,2}& \,=\, \dfrac{A}{\|\,U\,\|^{\,2}}\, \left\|\,U\,U^{\,\ast}\,f\,\right\|^{\,2} \,\leq\, A\, \left\|\,U^{\,\ast}\,f\,\right\|^{\,2}\\
& \,\leq\, \sum\limits_{\,j \,\in\, J}\, v_{j}^{\,2}\, \left\|\,\Lambda_{j}\,P_{\,W^{\,R}_{j}}\,\left(\,U^{\,\ast}\,f\,\right)\,\right\|^{\,2}\\
&\,=\, \sum\limits_{\,j \,\in\, J}\, v_{j}^{\,2}\, \left\|\,\Lambda_{j}\,P_{\,W^{\,R}_{j}}\,U^{\,\ast}\,P_{\,U\,W^{\,R}_{j}}\,(\,f\,)\,\right\|^{\,2}.
\end{align*} 
Therefore, \,$\Gamma$\; is a \,$U\,U^{\,\ast}$-$g$-fusion frame for \,$\mathbb{H}^{\,R}\,(\,\mathfrak{Q}\,)$. 
\end{proof}

\begin{theorem}
Let \,$U$\, be an bounded right $\mathfrak{Q}$-linear operator on \,$\mathbb{H}^{\,R}\,(\,\mathfrak{Q}\,)$\, and \,$\Gamma \,=\, \left\{\,\left(\,U\,W_{j},\, \Lambda_{j}\,P_{\,W_{j}}\,U^{\,\ast},\, v_{j}\,\right)\,\right\}_{j \,\in\, J}$\; be a \,$g$-fusion frame for \,$\mathbb{H}^{\,R}\,(\,\mathfrak{Q}\,)$.\,Then \,$\Lambda \,=\, \left\{\,\left(\,W_{j} ,\, \Lambda_{j},\, v_{j}\,\right)\,\right\}_{j \,\in\, J}$\, is a \,$g$-fusion frame for \,$\mathbb{H}^{\,R}\,(\,\mathfrak{Q}\,)$.   
\end{theorem}

\begin{proof}
For each \,$f \,\in\, \mathbb{H}^{\,R}\,(\,\mathfrak{Q}\,)$, we have
\begin{align*}
\dfrac{A}{\|\,U\,\|^{\,2}}\,\left \|\,f\,\right \|^{\,2}& \,=\, \dfrac{A}{\|\,U\,\|^{\,2}}\,\left\|\,U^{\,\ast}\,(\,U^{\,-\, 1}\,)^{\,\ast}\,f\,\right\|^{\,2} \,\leq\, A\;\left\|\,\left(\,U^{\,-\, 1}\,\right)^{\,\ast}\,f\,\right\|^{\,2}\\
& \leq\, \sum\limits_{\,j \,\in\, J}\, v_{j}^{\,2}\, \left\|\,\Lambda_{j}\,P_{\,W^{\,R}_{j}}\,U^{\,\ast}\,P_{\,U\,W^{\,R}_{j}}\,\left(\,\left(\,U^{\,-\, 1}\,\right)^{\,\ast}\,f\,\right)\,\right\|^{\,2}\\
& \,=\, \sum\limits_{\,j \,\in\, J}\, v_{j}^{\,2}\, \left\|\,\Lambda_{j}\,P_{\,W^{\,R}_{j}}\,\left(\,U^{\,\ast}\,\left(\,U^{\,-\, 1}\,\right)^{\,\ast}\,f\,\right)\,\right\|^{\,2}\\
& \,=\, \sum\limits_{\,j \,\in\, J}\, v_{j}^{\,2}\, \left\|\,\Lambda_{j}\,P_{\,W^{\,R}_{j}}\,(\,f\,)\,\right\|^{\,2}.
\end{align*}
Also, for each \,$f \,\in\, \mathbb{H}^{\,R}\,(\,\mathfrak{Q}\,)$, we have
\begin{align*}
\sum\limits_{\,j \,\in\, J}\, v_{j}^{\,2}\, \left\|\,\Lambda_{j}\,P_{\,W^{\,R}_{j}}\,(\,f\,)\,\right\|^{\,2}& \,=\, \sum\limits_{\,j \,\in\, J}\, v_{j}^{\,2}\, \left\|\,\Lambda_{j}\,P_{\,W^{\,R}_{j}}\,\left(\,U^{\,\ast}\,\left(\,U^{\,-\, 1}\,\right)^{\,\ast}\,f\,\right)\,\right\|^{\,2}\\
&\,=\, \sum\limits_{\,j \,\in\, J}\, v_{j}^{\,2}\, \left\|\,\Lambda_{j}\,P_{\,W^{\,R}_{j}}\,U^{\,\ast}\,P_{\,U\,W^{\,R}_{j}}\,\left(\,\left(\,U^{\,-\, 1}\,\right)^{\,\ast}\,f\,\right)\,\right\|^{\,2}\\
&\leq\, B \; \left\|\,\left(\,U^{\,-\, 1}\,\right)^{\,\ast}\,f\, \right\|^{\,2} \,\leq\, B\; \left\|\,U^{\,-\, 1}\,\right\|^{\,2}\,\|\,f\,\|^{\,2}.
\end{align*}
Thus, \,$\Lambda$\, is a \,$g$-fusion frame for \,$\mathbb{H}^{\,R}\,(\,\mathfrak{Q}\,)$\, with bounds \;$\dfrac{A}{\|\,U\,\|^{\,2}}$\, and \,$B\; \left\|\,U^{\,-\, 1}\,\right\|^{\,2}$.         
\end{proof}

\end{document}